\documentclass{amsart}
\usepackage{amssymb} 
\usepackage{amsmath} 
\usepackage{amscd}
\usepackage{amsbsy}
\usepackage{comment}
\usepackage[matrix,arrow]{xy}
\usepackage{hyperref}

\DeclareMathOperator{\Rad}{Rad}
\DeclareMathOperator{\Radp}{Rad_{\{2,5,7\}}}
\DeclareMathOperator{\Radpp}{Rad_{\{2,3,5,7\}}}
\DeclareMathOperator{\Norm}{Norm}

\DeclareMathOperator{\GL}{GL}

\DeclareMathOperator{\Res}{Res}

\DeclareMathOperator{\Gal}{Gal}

\DeclareMathOperator{\ord}{ord}

\DeclareMathOperator{\Aut}{Aut}

\newcommand{\Q}{{\mathbb Q}}
\newcommand{\Z}{{\mathbb Z}}
\newcommand{\F}{{\mathbb F}}

\newcommand{\cT}{\mathcal{T}}

\newcommand{\OO}{{\mathcal O}}
\newcommand{\ga}{{\mathfrak{a}}}
\newcommand{\gb}{{\mathfrak{b}}}

\newcommand{\sS}{\mathfrak{S}}
\newcommand{\fp}{\mathfrak{p}}
\newcommand{\fq}{\mathfrak{q}}
\newcommand{\fr}{\mathfrak{r}}
\newcommand{\fl}{\mathfrak{l}}
\newcommand{\mP}{\mathcal{P}}

\begin {document}

\newtheorem{thm}{Theorem}
\newtheorem{lem}{Lemma}[section]

\theoremstyle{definition}

\theoremstyle{remark}

\title[Superelliptic equations]{Superelliptic equations arising from sums of consecutive powers}
\author[Michael Bennett]{Michael A. Bennett}
\address{Department of Mathematics, University of British Columbia, Vancouver, B.C., V6T 1Z2 Canada}
\email{bennett@math.ubc.ca}

\author{Vandita Patel}
\address{Mathematics Institute, University of Warwick, Coventry CV4 7AL, United Kingdom}
\email{vandita.patel@warwick.ac.uk}

\author{Samir Siksek}
\address{Mathematics Institute, University of Warwick, Coventry CV4 7AL, United Kingdom}
\email{S.Siksek@warwick.ac.uk}
\thanks{
The first-named author is supported
by NSERC. 
The second-named author is supported by an EPSRC studentship.
The third-named author
is supported by 
the EPSRC {\em LMF: L-Functions and Modular Forms} Programme Grant
EP/K034383/1.
}

\date{\today}


\begin {abstract}
Using only elementary arguments, Cassels solved the 
Diophantine equation $(x-1)^3+x^3+(x+1)^3=z^2$ (with $x$, $z \in \Z$).
The generalization $(x-1)^k+x^k+(x+1)^k=z^n$
(with $x$, $z$, $n\in \Z$ and $n\ge 2$)
was considered by
Zhongfeng Zhang who solved it
for $k \in \{ 2, 3, 4 \}$ 
using Frey-Hellegouarch curves and their corresponding Galois representations.
In this paper, by employing some sophisticated refinements of this approach,
we show that the only solution for $k=5$ is $x=z=0$,
and that there are no solutions for $k=6$. The chief innovation we employ is  a computational one, which enables
us to avoid the full computation of data about cuspidal newforms of high level.
\end {abstract}
\maketitle

\section{Introduction} \label{intro}

In 1964, Leveque  \cite{Le} applied a theorem of Siegel \cite{Si} to show that, 
if $f(x) \in \mathbb{Z}[x]$ is a polynomial of degree $k
\geq 2$ with at least two simple roots, and $n \geq \max \{ 2, 5-k \}$
is an integer, then the {\it superelliptic}  equation
\begin{equation} \label{eq-1}
f(x) = z^n
\end{equation}
has at most finitely many solutions in integers $x$ and $z$.
This result was extended by  Schinzel and Tijdeman
\cite{ScTi}, through application of lower bounds for linear forms in logarithms, to show that equation (\ref{eq-1}) has in fact at
most finitely many solutions in integers $x, z$ and {\it variable}
$n \geq \max \{ 2, 5-k \}$ (where we count the solutions with $z^n = \pm 1,0$ once). 

While this latter result  is effective (in the sense that the finite set of
triples $(x,z,n)$ is effectively computable), in practice such a determination has  infrequently been achieved, due to the extraordinary size of the bounds for $x, z$ and $n$ arising from the proof. The few cases that have been treated in the literature have been restricted to polynomials with very few monomials, or with multiple linear factors over $\Q$.

One class of polynomials that has proved, in certain cases at least, amenable to such 
an approach, is that arising from sum of consecutive $k$-th powers.
Let us define 
$$
S_k(x) = 1^k + 2^k + \cdots + x^k,
$$
where $x$ and $k$ are nonnegative integers. Equations of the shape
\begin{equation} \label{super}
S_k(x) - S_k(y) = z^n
\end{equation}
have been considered by a number of authors, under the hypotheses that $y=0$ (see e.g. \cite{Br}, \cite{GTV}, \cite{Ha},
\cite{P}, \cite{P2}, \cite{Sc}, \cite{U1}, \cite{VGT}), that $y=[x/2]$ (\cite{Zhang}) and that $y=x-3$ (\cite{Cassels}, \cite{Zhang}). In the first two of these situations, the resulting polynomials on the left hand side of equation (\ref{super}) have at least two distinct linear factors over $\Q$, which allows the problem to be reduced to one of  binomial Thue equations.

Regarding the last of these cases, Cassels \cite{Cassels} solved the 
Diophantine equation $(x-1)^3+x^3+(x+1)^3=z^2$ in integers $x$ and $z$,
showing that the only solutions satisfy $x=0$, $1$, $2$ and $24$.
Zhongfeng Zhang \cite{Zhang} subsequently considered the more general
equation
\begin{equation}\label{eqn:Zhang}
(x-1)^k+x^k+(x+1)^k=z^n, \qquad x,~z,~k,~n\in \Z, \quad k,~n \ge 2.
\end{equation}
Associating solutions to a Frey-Hellegouarch curve and applying standard
level lowering arguments he proved that the only solutions
with $k \in \{ 2, 3, 4 \}$ are $(x,z,k,n)=(1,\pm 3, 3, 2)$,
$(2,\pm 6,3,2)$, $(24,\pm 204, 3,2)$, $(\pm 4, \pm 6, 3, 3)$ and
$(0,0,3,n)$. 

In this paper, we extend Zhang's result,  completely solving equation (\ref{eqn:Zhang}) in the cases $k =5$ and $k=6$. It should be emphasized that these results cannot apparently be obtained from the arguments of \cite{Zhang}, using Frey-Hellegouarch curves over $\Q$. Indeed, the purpose of this paper is two-fold. On the one hand, 
we will use the  case $k=5$ to advertise the utility of
the more powerful multi-Frey-Hellegouarch approach, first pioneered in \cite{multi}. Our result here is as follows :

\begin{thm}\label{thm:main}
The only solutions to the equation
\[
(x-1)^5+x^5+(x+1)^5=z^n, \qquad x,~z,~n \in \Z, \quad n \ge 1
\]
satisfy $x=z=0$.
\end{thm}

The other main purpose of this paper is to introduce a new computational approach to handle Diophantine problems where the problem of  extracting information about associated forms arising from modularity is at the limit of current computational power.
The method of Frey-Hellegouarch curves and Galois representations generally
requires the explicit computation of weight $2$ newforms of certain levels
and also the computation of some of their Hecke eigenvalues.
This computation can be completely impractical if the level is large. 
This turns out to be the case for equation  (\ref{eqn:Zhang}) with $k=6$, where required
newforms have level $3^3 \cdot 3391$ and the newform space has
dimension $4520$. We develop a version of the standard \lq method
for bounding exponents\rq\ \cite[Section 9]{Siksek} that does not
require the computations of the newforms, but merely a few
(computationally
much less expensive) Hecke 
polynomials. This allows us to prove the following theorem.
\begin{thm}\label{thm:2}
The equation
\[
(x-1)^6+x^6+(x+1)^6=z^n, \qquad x,~z,~n \in \Z, \quad n \ge 1
\]
has no solution.
\end{thm}

We are grateful to Steve Donnelly and John Voight for much help
and advice on the computation of newforms in \texttt{Magma}.

\section{The case $k=5$ and two Fermat equations of signature $(p,p,2)$}
The equation $(x-1)^5+x^5+(x+1)^5=z^n$ can be rewritten as
\begin{equation}\label{eqn:main}
x(3x^4+20x^2+10)=z^n \, .
\end{equation}
It suffices to deal with the case $n=p$ where $p$ is
a prime. 
We write $\alpha=\gcd(x,10)$, whereby
\begin{equation}\label{eqn:factor}
x=\alpha^{p-1} z_1^p \; \; \mbox{ and } \; \; 3x^4+20 x^2+10= \alpha z_2^p.
\end{equation}
We shall use this factorization to construct two associated
Fermat equations with signatures $(p,p,2)$.
We make use of the identity
\[
7 x^4+(3 x^4+20 x^2+10)=10 (x^2+1)^2.
\]
Substituting from \eqref{eqn:factor} and dividing by $\alpha$
we obtain
\begin{equation}\label{eqn:ferm1}
7 \alpha^{4p-5} z_1^{4p}+z_2^p=(10/\alpha) (x^2+1)^2.
\end{equation}
The reader will observe that this is a generalized Fermat
equation with signature $(p,p,2)$ where the three
terms are coprime. 

We also make use of the identity
\[
3(3 x^4+20 x^2+10)+70=(3x^2+10)^2.
\]
Again substituting from \eqref{eqn:factor} and dividing by $\alpha$
we obtain
\begin{equation}\label{eqn:ferm2}
3 z_2^p+\frac{70}{\alpha}=\alpha \left(\frac{3 x^2+10}{\alpha} \right)^2.
\end{equation}
Once again, the three terms in this equation are integral and coprime.
We interpret this as a generalized Fermat equation with signature
$(p,p,2)$ by treating the term $70/\alpha$ as $(70/\alpha) \cdot 1^p$.

We will associate a Frey-Hellegouarch curve to
each of the Fermat equations \eqref{eqn:ferm1} and 
\eqref{eqn:ferm2}, and use the information derived simultaneously
from
both Frey-Hellegouarch curves to prove Theorem~\ref{thm:main} for $n=p \ge 7$.
We need to treat exponents $p=2$, $3$ and $5$ separately;
 we do this in the next section.

\section{The case $k=5$: small values of $p$}

\begin{lem}
The only solutions to \eqref{eqn:main} with $n=p=2$ 
are $x=z=0$. 
\end{lem}
\begin{proof}
Write $X=3 \alpha x^2$ and $Y=3 \alpha x z_2$. From \eqref{eqn:factor},
it follows that $(X,Y)$ is an integral point on the elliptic curve
\[
E_\alpha \; : \; Y^2=X(X^2+20 \alpha X+30 \alpha^2).
\]
Using the computer algebra package \texttt{Magma} \cite{magma},
we determine the integral points on $E_\alpha$. For this computation,
\texttt{Magma}
applies the standard linear forms in elliptic logarithms method 
\cite[Chapter XIII]{Smart}. The integral points on
these curves are
$$
(-6 , \pm 18), \; \;  (-5 , \pm 15 ),  \; \;  (0 , 0 ) \; \mbox{ and } \; (1080 , \pm 35820 ),
$$
for $\alpha=1$; the points
$(-54 , \pm 306 )$ and $(0 , 0 )$
for $\alpha=5$; and just the point $(0,0)$ for $\alpha=2$ or $10$.
The lemma follows immediately.
\end{proof}

\begin{lem}
The only solutions to \eqref{eqn:main} with $n=p=3$ 
are $x=z=0$. 
\end{lem}
\begin{proof}
Let $X=3\alpha z_2$ and $Y=3\alpha (3 x^2+10)$. From \eqref{eqn:ferm2},
we see that $(X,Y)$ is an integral point on the elliptic curve
\[
E_\alpha \; : \; Y^2=X^3+630 \alpha^2 \, .
\]
Again using \texttt{Magma}, we determined the integral points on these four
elliptic curves. The curve $E_1$ has no integral points.
The integral points on $E_5$ are
$(-5 , \pm 125 )$ and $(99 , \pm 993 )$, while those on $E_2$
are 
$(-6 , \pm 48 ),  (9 , \pm 57 )$ and  $(46 , \pm 316 )$.
Finally, the integral points on $E_{10}$ are given by
$$
 (1 , \pm 251 ), \; \; 
(30 , \pm 300 ), \; \; 
(81 , \pm 771 ) \; \mbox{ and } \;  (330 , \pm 6000 ).
$$
The lemma follows.
\end{proof}

\begin{lem}\label{lem:peq5}
The only solutions to \eqref{eqn:main} with $n=p=5$ 
are $x=z=0$. 
\end{lem}
\begin{proof}
From \eqref{eqn:factor} we have
\[
x=\alpha^4 z_1^5, \qquad (3x^2+10+\sqrt{70})(3 x^2+10-\sqrt{70})=3 \alpha z_2^5.
\]
Let $K=\Q(\sqrt{70})$. This field has ring of integers $\OO_K=\Z[\sqrt{70}]$
and fundamental unit $\varepsilon=251+31 \sqrt{70}$. We consider the following prime ideals 
\begin{gather*}
\fp_2=(2,\sqrt{70}), \qquad \fp_3=(3,1+\sqrt{70}), \qquad \fp_3^\prime=(3,1-\sqrt{70}), \\
\fp_5=(25+3\sqrt{70})\OO_K \; \; \; \mbox{ and } \; \; \;  \fp_7=(7,\sqrt{70}).
\end{gather*}
These satisfy
\[
\fp_2^2=2\OO_K, \; \; \;  \fp_5^2=5\OO_K, \; \; \;  \fp_7^2=7\OO_K \; \mbox{ and } \;  \fp_3 \fp_3^\prime=3\OO_K.
\]
The field $K$ has class number $2$, with $\fp_2$, $\fp_3$, $\fp_3^\prime$ and $\fp_7$
all representing the non-trivial ideal class. Observe that
\[
\ord_{\fp_3}(10+\sqrt{70})=1 \; \; \mbox{ and } \; \;  \ord_{\fp_3}(10-\sqrt{70})=0.
\]
Moreover,
\[
\ord_{\fp_2}(10+\sqrt{70})=\ord_{\fp_5}(10+\sqrt{70})=1.
\]
Let $a=\ord_2(\alpha)$ and $b=\ord_5(\alpha)$, so that $a, b \in \{ 0, 1 \}$.  We deduce that
\[
(3x^2+10+\sqrt{70}) \OO_K=
\ga \cdot \gb^5, \; \; \mbox{ where } \; \;
\ga=\fp_2^a \cdot \fp_5^b \cdot \fp_3 
\]
and $\gb$ is an ideal of $\OO_K$. Observe that $\gb$ is principal if and only if $\ga$
is principal. 
Let
\[
\fq=\begin{cases}
1 \cdot \OO_K & \text{if $\ga$ is principal}\\
\fp_7 & \text{if $\ga$ is non-principal}.
\end{cases}
\]
Then we can write
\[
(3x^2+10+\sqrt{70}) \OO_K= (\ga \fq^{-5}) \cdot (\gb \fq)^5
\]
where both $\ga \fq^{-5}$ and $\gb \fq$ are principal; the former is a fractional ideal,
while the latter is an integral ideal. Write
\[
\ga \fq^{-5}=\frac{r+s \sqrt{70}}{d} \OO_K
\]
where $r$, $s$, $d \in \Z$, with $d \ge 1$ as small as possible. 
Now
\[
3x^2+10+\sqrt{70}=\frac{1}{d}(r+s\sqrt{70}) \cdot \varepsilon^{c} \cdot (u+v \sqrt{70})^5, \qquad -2 \le c \le 2
\]
with $u$, $v$ in $\Z$. Comparing coefficients of $1$, $\sqrt{70}$, and recalling that $x=\alpha^4 z_1^5$
we have
\begin{equation}\label{eqn:system}
f(u,v)=d(3 \alpha^8 z_1^{10}+10) \; \; \mbox{ and } \; \;  g(u,v)=d,
\end{equation}
where $f$, $g \in \Z[u,v]$ are homogenous of degree $5$. Observe that $d$ is determined by $\alpha$,
while $f$ and $g$ are determined by $\alpha$ and $c$. For each possibility for $\alpha$ and $c$
we checked the system  \eqref{eqn:system} for solubility modulo $2^6$, $3^3$, $5^3$, $7^3$
 and all primes $11 \le q < 100$. This allowed us to eliminate all possibilities except for 
$(\alpha,c)=(2,2)$ and $(\alpha,c)=(10,0)$. For both these possibilies $d=1$. The second equation
in \eqref{eqn:system} is in fact a Thue equation. 
We used \texttt{Magma} to solve both Thue equations; for the theory behind \texttt{Magma}'s Thue
equation solver see \cite[Chapter VII]{Smart}.

For $(\alpha,c)=(2,2)$ this Thue equation is
\[
5521 u^5 + 230960 u^4 v + 3864700 u^3 v^2 + 32334400 u^2 v^3 + 135264500 u v^4 + 226340800 v^5=1,
\]
and we found that it has no solutions. 
For $(\alpha,c)=(10,0)$, the corresponding Thue equation is
\[
u^5 + 50 u^4 v + 700 u^3 v^2 + 7000 u^2 v^3 + 24500 u v^4 + 49000 v^5=1.
\]
The only solution is $(u,v)=(1,0)$. Since the first equation in \eqref{eqn:system} is, in this case,
\[
10 u^5 + 350 u^4 v + 7000 u^3 v^2 + 49000 u^2 v^3 + 245000 u v^4 + 343000v^5=3 \cdot 10^8 \cdot z_1^{10}+10,
\]
it follows that $z_1=0$ and hence $x=0$ as required.
\end{proof}

\section{The case $k=5$: first Frey-Hellegouarch curve}
Henceforth we suppose that $p \ge 7$ and
that $x \ne 0$.
We apply the recipes of the first author and Skinner
\cite[Section 2]{BenS} to equation \eqref{eqn:ferm1}
(see also \cite{Siksek}; this latter reference
is a comprehensive tutorial on the modular approach).
The recipes lead us to
attach to \eqref{eqn:ferm1} a Frey-Hellegouarch curve $E_{x,\alpha}$
which depends on $\alpha$ as well as $x$. The possible
values for $\alpha$ are $1$, $5$, $2$ and $10$. The 
corresponding Frey-Hellegouarch elliptic curves are
\begin{gather}
\label{eqn:Frey-Hellegouarch1}
E_{x,1} \; : \; 
Y^2=X^3+20 (x^2+1)X^2 + 10(3 x^4+20 x^2+10) X  \\
\label{eqn:Frey-Hellegouarch2}
E_{x,5} \; : \; 
Y^2=X^3+4 (x^2+1) X^2 + \frac{2(3 x^4+20 x^2+10)}{5} X \\
\label{eqn:Frey-Hellegouarch3}
E_{x,2} \; : \; 
Y^2+XY=X^3 +\frac{(5 x^2+4)}{4} X^2 + \frac{35 x^4}{128} X \\
\label{eqn:Frey-Hellegouarch4}
E_{x,10} \; : \; 
Y^2+XY=X^3+\frac{x^2}{4} X^2 + \frac{7 x^4}{640} X \, .
\end{gather} 
For a non-zero integer $u$ and a set of primes $S$, we define $\Rad_S(u)$ to be
the product of the distinct prime divisors of $u$ that do not belong to $S$.
For an elliptic curve $E/\Q$, we denote
its minimal discriminant and conductor by $\Delta(E)$
and $N(E)$.
\begin{lem}
The elliptic curves $E_{x,\alpha}$ have non-trivial $2$-torsion
over $\Q$.
Their discriminant and conductors are  
$$
\begin{array}{l}
\Delta(E_{x,1})=2^9 \cdot 5^3 \cdot 7 \cdot z_1^{4p} \cdot z_2^{2p},
\; \; \; 
N(E_{x,1})=2^8 \cdot 5^2 \cdot 7 \cdot \Radp(z_1 z_2)\\
\Delta(E_{x,5})=2^9 \cdot 5^{4p-5} \cdot 7 \cdot z_1^{4p}\cdot z_2^{2p},
\; \; \; 
N(E_{x,5})=2^8 \cdot 5 \cdot 7 \cdot \Radp(z_1 z_2)\\
\Delta(E_{x,2})=2^{8p-22} \cdot 5^3 \cdot 7^2 \cdot z_1^{8p}\cdot z_2^{p},
\; \; \; 
N(E_{x,2})=2 \cdot 5^2 \cdot 7 \cdot \Radp(z_1 z_2)\\
\Delta(E_{x,10})=2^{8p-22} \cdot 5^{8p-10} \cdot 7^2 z_1^{8p} \cdot z_2^{p},
\; \; \; 
N(E_{x,10})=2 \cdot 5 \cdot 7 \cdot \Radp(z_1 z_2) \, .
\end{array}
$$
\end{lem}
\begin{proof}
This follows from \cite[Lemma 2.1]{BenS}.
\end{proof}
We note in passing that we have already used the assumption $x \ne 0$.
If $x=0$, then $z_1=0$ and the curves $E_{x,\alpha}$ are not
elliptic curves but merely singular Weierstrass equations
(i.e. the discriminant $\Delta(E_{x,\alpha})=0$). We maintain
the assumption $x \ne 0$ throughout. 

For an elliptic curve $E/\Q$,
we write $\overline{\rho}_{E,p}$ for the mod $p$ representation
giving the action of $G_\Q=\Gal(\overline{\Q}/\Q)$ on the $p$-torsion
$E[p]$:
\[
\overline{\rho}_{E,p} \;
: \; 
G_\Q \rightarrow \Aut(E[p]) \cong \GL_2(\F_p).
\]
If $\overline{\rho}_{E,p}$ arises from a newform $f$
then we write $E \sim_p f$.
\begin{lem}\label{lem:ll1}
Let $E_{x,\alpha}$ be the Frey-Hellegouarch curve in \eqref{eqn:Frey-Hellegouarch1}--\eqref{eqn:Frey-Hellegouarch4}.
Then $E_{x,\alpha}\sim_p f$ where $f$ is a newform
of weight $2$ and level $L_\alpha$:
\[
L_1=2^8 \cdot 5^2 \cdot 7,
\; \; \; 
L_5=2^8 \cdot 5 \cdot 7,
\; \; \; 
L_2=2 \cdot 5^2 \cdot 7
\; \; \mbox{ and } \; \; 
L_{10}=2 \cdot 5 \cdot 7.
\]
\end{lem}
\begin{proof}
This is immediate from \cite[Lemma 3.2]{BenS} which
in turn relies on modularity of elliptic curves over $\Q$
due to Wiles, Breuil, Conrad, Diamond and Taylor \cite{Wiles}, \cite{Mod},
on Ribet's level lowering theorem \cite{Ribet} and also on
irreducibility theorems for mod $p$ representations of elliptic
curves due to Mazur \cite{Mazur}. It is here that the
assumption $p \ge 7$ is used to ensure the irreducibility
of the representation $\overline{\rho}_{E_{x,\alpha},p}$.
\end{proof}
Using \texttt{Magma}, we computed the weight $2$ newforms
of levels $L_\alpha$. The results of this computation 
are summarized in 
Table~\ref{table:ll1}. For this computation we used
\texttt{Magma}'s highly optimized
Hilbert modular forms package (the classical
newforms we are computing can be regarded as Hilbert newforms over $\Q$).
The theory and algorithms behind this package are described in
\cite{DV}.

\begin{table}
\begin{centering}
\begin{tabular}{|c|c|c|c|}
\hline\hline
$\alpha$ &  $\dim S_2^{\mathrm{new}}(L_\alpha)$ & number of conjugacy & $(d, \text{number of newforms
of degree $d$})$\\
         &                                      & classes of forms &\\
\hline\hline
$1$      & $912$ & $196$ & $(1, 52)$, $(2, 32)$, $(3, 12)$, $(4, 22)$, \\
& & & $(5, 8)$, $(6, 28)$, $(8, 12)$, $(9, 8)$, \\
& & & $(12, 16)$, $(16, 2)$, $(18, 4)$\\
\hline
$5$      & $192$  & $64$ &  $(1, 20)$, $(2, 12)$, $(3, 12)$, $(4, 4)$,\\ 
& & & $(6, 16)$\\
\hline
$2$   & $10$ & $8$ & $(1, 6)$, $(2, 2)$ \\
\hline
$10$  & $1$  & $1$ & $(1,1)$ \\
\hline
\end{tabular}
\caption{Information for weight $2$ newforms of level
$L_\alpha$, where $L_\alpha$ is given by Lemma~\ref{lem:ll1}.}
\label{table:ll1}
\end{centering}
\end{table}

\section{The case $k=5$: second Frey-Hellegouarch curve}
Applying the recipes of Bennett and Skinner \cite[Section 2]{BenS} to equation \eqref{eqn:ferm2}
leads
us to associate to this the Frey-Hellegouarch elliptic curve
\begin{equation}\label{eqn:SecondFrey-Hellegouarch}
F_{x,\alpha} \; : \; Y^2=X^3+2(3x^2+10)X^2+70 X \, .
\end{equation}
Although the equation for this Frey-Hellegouarch curve is independent of $\alpha$,
the discriminant and conductor do depend on $\alpha$.
\begin{lem}
The elliptic curve $F_{x,\alpha}$ has non-trivial $2$-torsion
over $\Q$.
Its discriminant and conductor are given by
\begin{gather*}
\Delta(F_{x,1})=2^8 \cdot 3 \cdot 5^2 \cdot 7^2 \cdot z_2^{p},
\qquad
N(F_{x,1})=2^7 \cdot 3 \cdot 5 \cdot 7 \cdot \Radpp(z_2)\\
\Delta(F_{x,5})=2^8 \cdot 3 \cdot 5^3 \cdot 7^2 \cdot z_2^{p},
\qquad
N(F_{x,5})=2^7 \cdot 3 \cdot 5^2 \cdot 7 \cdot \Radpp(z_2)\\
\Delta(F_{x,2})=2^9 \cdot 3 \cdot 5^2 \cdot 7^2 \cdot z_2^{p},
\qquad
N(F_{x,2})=2^8 \cdot 3 \cdot 5 \cdot 7 \cdot \Radpp(z_2)\\
\Delta(F_{x,10})=2^9 \cdot 3 \cdot 5^3 \cdot 7^2 z_2^{p},
\qquad
N(F_{x,10})=2^8 \cdot 3 \cdot 5^2 \cdot 7 \cdot \Radpp(z_2) \, .
\end{gather*}
\end{lem}
\begin{proof}
Again this follows from \cite[Lemma 2.1]{BenS}.
\end{proof}

\begin{lem}\label{lem:ll2}
Let $F_{x,\alpha}$ be the Frey-Hellegouarch curve in \eqref{eqn:SecondFrey-Hellegouarch}.
Then $F_{x,\alpha}\sim_p g$ where $g$ is a newform
of weight $2$ and level $M_\alpha$, where
\[
M_1=2^7 \cdot 3 \cdot 5 \cdot 7,
\; \; 
M_5=2^7 \cdot 3 \cdot 5^2 \cdot 7,
\; \;
M_2=2^8 \cdot 3 \cdot 5 \cdot 7
\; \mbox{ and } \; 
M_{10}=2^8 \cdot 3 \cdot 5^2 \cdot 7.
\]
\end{lem}
\begin{proof}
Again this is immediate from \cite[Lemma 3.2]{BenS}.
\end{proof}
Table~\ref{table:ll2} gives information about the spaces
of newforms of weight $2$ and level $M_\alpha$. 

\begin{table}
\begin{centering}
\begin{tabular}{|c|c|c|c|}
\hline\hline
$\alpha$ &  $\dim S_2^{\mathrm{new}}(M_\alpha)$ & number of conjugacy & $(d, \text{number of newforms
of degree $d$})$\\
         &                                      & classes of forms &\\
\hline\hline
$1$      & $192$ & $112$ & 
$(1, 64)$, $(2, 28)$, $(3, 12)$, $(4, 4)$,\\
& & &  $(5, 4)$\\
\hline
$5$      & $912$ & $356$ &  $(1, 176)$, $(2, 64)$, $(3, 12)$, $(4, 36)$, \\
& & & $(5, 28)$, $(6, 8)$, $(7, 24)$, $(9, 8)$\\
\hline
$2$    & $384$ & $128$ & $(1, 48)$, $(2, 16)$, $(3, 16)$, $(4, 28)$, \\
& & & $(6, 8)$, $(8, 12)$\\
\hline
$10$ & $1824$ & $396$ & $(1, 124)$, $(2, 60)$, $(3, 20)$, $(4, 52)$, \\
& & & $(5, 8)$, $(6, 40)$, $(8, 28)$\\
\hline
\end{tabular}
\caption{Information for weight $2$ newforms of level
$M_\alpha$, where $M_\alpha$ is given by Lemma~\ref{lem:ll2}.}
\label{table:ll2}
\end{centering}
\end{table}

\section{Proof of Theorem~\ref{thm:main}}
The following standard lemma \cite[Proposition 5.1]{Siksek} 
will be helpful in exploiting Lemmata~\ref{lem:ll1} and
\ref{lem:ll2}.
\begin{lem}\label{lem:info}
Let $E/\Q$ be an elliptic curve of conductor $N$ and $f=q+\sum_{i \ge 2} c_i q^i$
be a newform of weight $2$ and level $N^\prime \mid N$.
Write $K=\Q(c_1,c_2,\dots)$ for the totally real number field
generated by the Fourier coefficients of $f$. Suppose $E \sim_p f$
for some prime $p$. Then there is some prime ideal $\fp \mid p$
of $K$ such that for all primes $\ell$,
\begin{itemize}
\item if $\ell \nmid pN N^\prime$ then $a_\ell(E) \equiv c_\ell \pmod{\fp}$,
\item if $\ell \nmid p N^\prime$ and $\ell \mid \mid N$ then
$\pm (\ell+1) \equiv c_\ell \pmod{\fp}$.
\end{itemize}
\end{lem}

\bigskip

Fix a possible value for $\alpha \in \{ 1, 2, 5, 10 \}$. 
For convenience, we write $E_x$ and $F_x$ 
for the curves $E_{x,\alpha}$ and $F_{x,\alpha}$.
Note that the levels
$L_\alpha$ and $M_\alpha$ in Lemmata~\ref{lem:ll1} and~\ref{lem:ll2}
depend only on $\alpha$. 
Now fix a weight $2$ newform 
$f=q+\sum c_i q^i$ of level
$L_\alpha$ and another $g=q+\sum d_i q^i$ of level $M_\alpha$.
Suppose $E_{x} \sim_p f$ and $F_{x} \sim_p g$.
Write $K_1=\Q(c_1,c_2,\dotsc)$ and 
$K_2=\Q(d_1,d_2,\dotsc)$, and let $\ell > 7$ be a prime.
We would like to apply Lemma~\ref{lem:info} 
to obtain information about $p$. Suppose for now that $\ell \ne p$.
The Frey-Hellegouarch curves $E_{x}$ and $F_{x}$ depend
on the unknown $x$. However, their traces modulo $\ell$ 
depend only on $x$ modulo $\ell$.
Let $0 \le a \le \ell-1$ and suppose $x \equiv a \pmod{\ell}$.
We shall write $\Delta_1(x)$ for the discriminant of the
Weierstrass model $E_{x}$ and $\Delta_2(x)$ for the discriminant
of the Weierstrass model $F_{x}$ (these are polynomials in $x$).
Let
\begin{equation}\label{eqn:R}
R_\ell(f,a)=
\begin{cases}
\Norm_{K_1/\Q}((\ell+1)^2-c_\ell^2) & \text{if $\ell \mid \Delta_1(a)$}\\
\Norm_{K_1/\Q}(a_\ell(E_a)-c_\ell) & \text{if $\ell \nmid \Delta_1(a)$}.
\end{cases}
\end{equation}
It follows from Lemma~\ref{lem:info} and Lemma~\ref{lem:ll1} that
$p \mid R_\ell(f,a)$.
Let
\[
S_\ell(g,a)=
\begin{cases}
\Norm_{K_2/\Q}((\ell+1)^2-d_\ell^2) & \text{if $\ell \mid \Delta_2(a)$}\\
\Norm_{K_2/\Q}(a_\ell(F_a)-c_\ell) & \text{if $\ell \nmid \Delta_2(a)$}.
\end{cases}
\]
It further follows from Lemma~\ref{lem:info} and Lemma~\ref{lem:ll2} that
$p \mid S_\ell(g,a)$.
Now let
\[
T_\ell(f,g,a)=\gcd(R_\ell(f,a),S_\ell(g,a)).
\]
Then $p \mid T_\ell(f,g,a)$. Observe that while $a$ is unknown, as it is the residue
of $x$ modulo $\ell$, we may suppose that $0 \le a \le \ell -1$. Let
\[
T_\ell(f,g)=\ell \prod_{0 \le a \le \ell -1} T_\ell(f,g,a).
\]
Then $p \mid T_\ell(f,g)$. We had assumed above that $\ell \ne p$.
However as $\ell$ is a factor in the product defining $T_\ell(f,g)$,
the conclusion $p \mid T_\ell(f,g)$ is true even if $\ell=p$.
Finally we let
\[
U(f,g)=\gcd_{11 \le \ell < 100} T_\ell(f,g)
\]
where the gcd is taken over all primes $\ell$ in the range $11 \le \ell < 100$.
It follows that $p \mid U(f,g)$. To complete the proof of Theorem~\ref{thm:main},
we employ a simple \texttt{Magma} script that computes for each pair $(f,g)$
the quantity $U(f,g)$ and verifies that it is not divisible by
primes $\ge 7$. The computation took roughly $4$ days 
on a 2500MHz AMD Opteron, dominated by
the computation of the newforms.

\bigskip

\noindent \textbf{Remark.}
It is appropriate to comment at this stage as to whether the single Frey-Hellegouarch approach
(using either of the Frey-Hellegouarch curves $E_{x,\alpha}$ or $F_{x,\alpha}$
on its own) would have allowed us to establish Theorem~\ref{thm:main}.
The above argument is a multi-Frey-Hellegouarch version of the standard
single Frey-Hellegouarch method for bounding exponents (see \cite[Section 9]{Siksek}).
With notation as above, let
\[
B_\ell(f)= \ell \prod_{0 \le a \le \ell-1} R_\ell(f,a),
\]
for $\ell \ne 2$, $5$, $7$ (note that $3$ does not divide the possible
levels of $f$).
Under the assumption $E_x \sim f$,
the single Frey-Hellegouarch method for bounding exponents asserts that $p \mid B_\ell(f)$
and succeeds in bounding $p$ 
if we can find a prime $\ell \ne 2$, $5$, $7$ such that $B_\ell(f) \ne 0$.
Likewise, let
\[
C_\ell(g)=\ell \prod_{0 \le a \le \ell-1} S_\ell(g,a), 
\]
for $\ell \ne 2$, $3$, $5$, $7$.
Under the assumption $F_x \sim g$,
we have $p \mid C_\ell(g)$.
We first note that the solution $(x,z,n)=(0,0,p)$ of equation~\eqref{eqn:main}
leads to an elliptic curve  $F_{0,10}$ (i.e. a non-singular Weierstrass 
equation)
with Cremona reference
\texttt{134400BG1}. 
 Let $g$ be the eigenform (of level $M_{10}=134400$)
corresponding to the elliptic curve $F_{0,10}$.
Then $a_\ell(F_{0,10})=d_\ell$ where $g=q+\sum d_i q^i$.
Hence $S_\ell(g,0)=0$ and so $C_\ell(g)=0$ for all possible $\ell$.
Thus the single Frey-Hellegouarch method with the second Frey-Hellegouarch curve $F_{x,\alpha}$
fails to bound the exponent $p$.

The single Frey-Hellegouarch approach succeeds with the first Frey-Hellegouarch curve $E_{x,\alpha}$
in the sense that for all possible eigenforms $f$, we are able to find
some prime $\ell \ne 2$, $5$, $7$ such that $B_\ell(f) \ne 0$.
For any $\ell$, the bound $B_\ell(f)$ can be very large (especially if the
field of coefficients of $f$ has large degree). However, we consider 
instead
\[
B(f)=\gcd_{\ell \in \{3,11,13,\dotsc,97\} } B_\ell(f).
\]
If $E_x \sim f$ then $p \mid B(f)$. We computed the $B(f)$ for the possible
newforms $f$ and found many of them to be divisible by $7$ and $13$
though not by larger primes. It is possible to reduce the
cases $p=7$ and $p=13$ to Thue equations as in the proofs of 
Lemma~\ref{lem:peq5}. However the coefficients of these
Thue equations will be so unpleasant that we
do not expect to be able to solve them (uncondionally). 

\section{Dealing with small exponents for the case $k=6$}
We now consider the equation
\[
(x-1)^6+x^6+(x+1)^6=z^n, \qquad x,z,n \in \Z, \qquad n \ge 2,
\]
which corresponds to the case $k=6$ of \eqref{eqn:Zhang}.
This can be rewritten as 
$$
3x^6+30 x^4+30 x^2+2=z^n,
$$
whereby necessarily
$z^n \equiv 2 \pmod{3}$ and hence $n$ is odd. 
Moreover the polynomial $3t^6+30 t^4+30 t^2+2$ only takes
values $2$ and $3$ as $t$ ranges over $\F_7$. As these values
are not cubes in $\F_7$, we see that $3 \nmid n$. Thus
to prove Theorem~\ref{thm:2} for $k=6$ it is sufficient to
show that the equation
\begin{equation}\label{eqn:k6}
3x^6+30 x^4+30 x^2+2=z^p
\end{equation}
has no solutions with prime exponent $p \ge 5$.

\begin{lem}\label{lem:k6}
Equation \eqref{eqn:k6} has no solutions with $p=5$, $7$, $11$, $13$.
\end{lem}
\begin{proof}
Write $f=3t^6+30 t^4+30 t^2+2$. The polynomial $f$ is irreducible
over $\Q$. Let $\theta$ be a root of the equation $f(t)=0$, write $K=\Q(\theta)$
and let $\OO_K$ be the ring integers of $K$. 
The field $K$ has unit rank $2$ with $-1$ as a generator
for the roots of unity, and class group  $\cong (\Z/2\Z)^3 \times (\Z/36\Z)$.

Let $g(t)=f(t)/(x-\theta) \in K[t]$.
There are prime ideals $\fp$, $\fq$, $\fr_1$, $\fr_2$, $\fr_3$, $\fr_4$
such that
\begin{gather*}
2\cdot \OO_K=\fp^6, \qquad 3\cdot \OO_K=\fq^6,
\qquad
3391 \cdot \OO_K=\fr_1^2 \fr_2^2 \fr_3 \fr_4,\\
\theta \cdot \OO_K=\fp \fq^{-1} \; \; \; \mbox{ and } \; \; \; 
g(\theta) \cdot \OO_K=\fp^{11} \fq^3 \fr_1 \fr_2.
\end{gather*}

From \eqref{eqn:k6}, we know that
\[
(x-\theta) g(x)=z^p.
\]
Now $\ord_\fq(\theta)=-1$. As $x \in \Z$, we have $\ord_\fq(x-\theta)=-1$.
Let $\mP \ne \fq$ be a prime ideal and suppose that $\ord_\mP(x-\theta) 
\not \equiv 0 \pmod{p}$, whence $\ord_\mP(g(x)) \not\equiv 0 \pmod{p}$.
From the factorisation of $\theta \cdot \OO_K$
we know that $\ord_\mP(x-\theta)>0$. 
It is easy to see that $\mP \mid g(\theta)$. But $\ord_2(z)=0$
so $\mP=\fr_1$ or $\fr_2$.
Let $S=\{\fq,\fr_1,\fr_2\}$.
Hence $(x-\theta) {K^*}^p$ belongs to the \lq $p$-Selmer group\rq\
\[
K(S,p)=\{\alpha \in K^*/{K^*}^p \; : \; 
\text{$\ord_\mP(\alpha) \equiv 0 \pmod{p}$ for all $\mP \notin S$}
\}.
\]
This is an $\F_p$-vector space of finite dimension and,
for a given $p$, easy
to compute from the class group and unit group information (see 
\cite[Proof of Proposition VIII.1.6]{Silverman}).
Let
\[
\sS_p=\{ \alpha \in K(S,p) \; : \; \Norm(\alpha)=(1/3) {\Q^*}^p \}.
\]
Observe that $\Norm(x-\theta)=z^p/3$ so that $x-\theta \in \sS_p$.
Using \texttt{Magma}, we compute $K(S,p)$ and $\sS_p$ for $p=5$, $7$, $11$, 
$13$.
In all cases, $K(S,p)$ has $\F_p$-dimension equal to $5$, and the
set $\sS_p$ has $p^3$ elements. 

It follows that $x-\theta=\alpha \xi^p$ for some 
$\alpha \in \sS_p$ and $\xi \in K^*$.
We are now in a position to finally obtain a contradiction. 
Fix an $\alpha \in \sS_p$ such that $x-\theta=\alpha \xi^p$.
Let
$\ell \ne 3$ be a rational prime and $\fl_1,\dotsc,\fl_r$
be the prime ideals of $K$ dividing it. Suppose that none of the
$\fl_i$ belong to the support of $\alpha$. Let $x \equiv a \pmod{\ell}$
where $a \in \{0,1,\dotsc,\ell-1\}$. Then 
$(a-\theta)/\alpha \equiv \xi^p \pmod{\fl_i}$ for $i=1,\dotsc,r$.
Thus we may eliminate $\alpha$ if for each $a \in \{0,1,\dots,\ell-1\}$
there is some $i$ such $(a-\theta)/\alpha$ is not a $p$-th power
modulo $\fl_i$. For this to succeed,  
$\#\F_{\ell_i}=\Norm(\fl_i)$ needs to be $\equiv 1 \pmod{\ell}$.
For $p=5$, $7$, $11$, $13$ we have been able to find
a set of primes $\ell$ that we denote by $\cT_p$ which
allow us to eliminate all $\alpha \in \sS_p$. The sets $\cT_p$
are recorded in Table~\ref{table:k6}.

\begin{table}
\begin{centering}
\begin{tabular}{|l|c|}
\hline
$p$ & $\cT_p$ \\
\hline\hline
$5$ & $\{ 11, 191, 251, 691\}$\\
\hline
$7$ & $\{11, 337, 421, 491, 547\}$\\
\hline
$11$ & $\{397, 727, 859\}$\\
\hline
$13$ & $\{859, 1249\}$\\
\hline
\end{tabular}
\caption{The sets $\cT_p$ appearing in the proof of Lemma~\ref{lem:k6}}
\label{table:k6}
\end{centering}
\end{table}
\end{proof}

\section{Frey-Hellegouarch Curve for case $k=6$}

In this section, we construct a Frey-Hellegouarch curve attached to the equation~\eqref{eqn:k6}. 
In view of the previous section, we may suppose that the exponent $p$
in \eqref{eqn:k6} is a prime $\ge 17$.
The first author and Dahmen \cite{BD} attach a Frey-Hellegouarch curve to any equation of the form
$F(u,v)=z^p$ where $F$ is a  homogenous cubic form.
We now reproduce their recipe. Let
\[
H(u,v)=-\frac{1}{4}
\begin{vmatrix}
F_{uu} & F_{uv}\\
F_{uv} & F_{vv}
\end{vmatrix}  \; \; \mbox{ and } \; \; 
G(u,v)=\begin{vmatrix}
F_u & F_v \\
H_u & H_v
\end{vmatrix}.
\]
Associate to the solution $(u,v,z)$ of the equation $F(u,v)=z^p$
the Frey-Hellegouarch elliptic curve
\begin{equation}\label{eqn:BD}
E_{u,v}^\prime \; : \; Y^2=X^3-3 H(u,v) X +G(u,v) \, .
\end{equation}
This model has discriminant $2^4 \cdot 3^6 \cdot \Delta_F \cdot z^{2p}$,
where $\Delta_F$ is the discriminant of the binary form $F$.
Now consider the homogenous cubic form
\[
F(u,v)=3 u^3+30 u^2 v+ 30 u v^2+2 v^3.
\]
We note that $F(x^2,1)=3x^6+30 x^4+30 x^2+2$. 
Thus we may obtain a Frey-Hellegouarch curve for \eqref{eqn:k6}
by letting $(u,v)=(x^2,1)$ in \eqref{eqn:BD}.
In turns out that the model $E_{x^2,1}^\prime$
has bad reduction at $2$, but its quadratic twist by $2$
has good reduction at $2$, and we choose this to be the 
Frey-Hellegouarch curve associated to \eqref{eqn:k6}. 
A model which is
minimal at $2$ for this Frey-Hellegouarch curve is
\begin{multline}\label{eqn:Frey-Hellegouarch6}
E_x \; : \; Y^2 + Y \\ = X^3 + \frac{(-945 x^4 - 1269 x^2 - 1080)}{2}X +
\frac{(-15093 x^6 - 18630 x^4 + 26730 x^2 + 19061)}{4}.
\end{multline}
\begin{lem}\label{lem:minmod}
The model $E_x$ is integral, minimal and has discriminant and conductor
\[
\Delta_x=3^9 \cdot 3391 \cdot z^{2p}
\; \; \mbox{ and } \; \; N=3^3\cdot 3391 \cdot \Rad_{\{3,3391\}}(z) \, .
\]
\end{lem}
\begin{proof}
It is clear from \eqref{eqn:k6} that $x$ is odd,
whence one deduces that the model
$E_x$ is integral.
The discriminant for this model is
\[
\Delta_x=3^9 \cdot 3391 \cdot (3 x^6 + 30 x^4 + 30 x^2 + 2)^2=3^9 \cdot 3391 \cdot z^{2p},
\]
and the usual $c_4$-invariant is 
\[
c_4=2^3 \cdot 3^4 \cdot (35 x^4 + 47 x^2 + 40).
\]
We find that
\[
\Res(c_4,\Delta_x)=2^{40}\cdot 3^{84} \cdot 3391^{12} \; .
\]
Thus $E_x$ is minimal and semistable except possibly at
$p \in \{ 2, 3, 3391 \}$.
Since $\Delta_x$ is odd, $E_x$ in fact has good reduction at $2$.
We now show that $E_x$ has multiplicative reduction at $3391$.
The solutions to $c_4 \equiv 0 \pmod{3391}$ are
\[
x \equiv 983, \; 2408 \pmod{3391}.
\]
Both of these are roots to $3 x^6 + 30 x^4 + 30 x^2 + 2$
modulo $3391$. However for a solution $(x,z)$ to \eqref{eqn:k6},
we know that 
\[
3 x^6 + 30 x^4 + 30 x^2 + 2 \equiv 0 \pmod{3391^2}.
\]
We checked that $983$, $2408$ do not lift to roots for this
congruence. Hence $3391 \nmid c_4$. It follows that $E_x$
has multiplicative reduction at $3391$. Applying Tate's algorithm
\cite[Chapter IV]{SilvermanII},
we found that the $E_x$ has reduction type $\mathrm{IV}^*$ at $3$
with the valuation of the conductor equal to $3$.
The lemma follows.
\end{proof}


\begin{lem}\label{lem:irred}
Let $(x,z,p)$ be a solution to \eqref{eqn:k6} with $p \ge 17$
prime.
Let $E=E_x$ as in \eqref{eqn:Frey-Hellegouarch6}. 
Then $\overline{\rho}_{E,p}$
is irreducible.
\end{lem}
\begin{proof}
%
Suppose $\overline{\rho}_{E,p}$ is reducible. As $p \ge 17$,
it follows from the proof of Mazur's famous theorem on isogenies of elliptic
curves that
the $j$-invariant of $E$ belongs to $\Z[1/2]$ (see \cite[Corollary 4.4]{Mazur}).
However, $E$ has multiplicative reduction at $3391$ and so
$3391$ appears in the denominator of its $j$-invariant.
This contradiction shows that $\overline{\rho}_{E,p}$
is irreducible. 
\end{proof}

\section{Proof of Theorem~\ref{thm:2}}

\begin{lem}\label{lem:ll}
Let $(x,z,p)$ be a solution to \eqref{eqn:k6} with prime exponent
$p \ge 17$. Then
$E_x \sim_p f$ for some newform $f$ of weight $2$
and level $3^3 \cdot 3391$.
\end{lem}
\begin{proof}
This follows from Lemmata~\ref{lem:minmod} and~\ref{lem:irred}
together with Ribet's level lowering theorem \cite{Ribet}
(the special case \cite[Section 5]{Siksek} is enough for our
purpose).
\end{proof}
From Cremona's Database \cite{Cremona}, there are precisely four
elliptic curves having conductor $3^3 \cdot 3391$:
$$
\begin{array}{l}
F_1 \; : \;  y^2 + y = x^3 + 405 x + 22673\\ 
F_2 \; : \;  y^2 + y = x^3 + 45 x - 840\\ 
F_3 \; : \; y^2 + y = x^3 - 42 x - 104\\ 
F_4 \; : \;  y^2 + y = x^3 - 378 x + 2801 .
\end{array}
$$
\begin{lem}\label{lem:notrat}
$E_x \not \sim_p F_i$ for 
$i=1$, $2$, $3$, $4$.
\end{lem}
\begin{proof}
Suppose $E_x \sim F_i$.
As $2$ is a prime of good reduction for both elliptic curves,
Then $a_2(E_x) \equiv a_2(F_i) \pmod{p}$.
From \eqref{eqn:Frey-Hellegouarch6} and the fact that $x$ is odd, we find that
\[
E/\mathbb{F}_2 \; : \; Y^2+Y=X^3+X+1 \; .
\]
It follows that $a_2(E_x)=2$.  Since
\[
a_2(F_1)=2, \; \;  a_2(F_2)=-2, \mbox{ and }  a_2(F_3)= a_2(F_4)=0,
\]
we thus have $i=1$. 

Next we apply the method
of bounding the exponents.  For a prime $\ell \ne 3$, $3391$, let
\begin{equation}\label{eqn:R1}
R_\ell(a)=
\begin{cases}
(\ell+1)^2-a_\ell(F_1)^2 & \text{if $\ell \mid \Delta_a$}\\
a_\ell(E_a)-a_\ell(F_1) & \text{if $\ell \nmid \Delta_a$},
\end{cases}
\end{equation}
and
\[
B_\ell= \ell \prod_{0 \le a \le \ell-1} R_\ell(a).
\]
It follows from Lemma~\ref{lem:info} that $p \mid B_\ell$.
We find that
$B_{11}=5^4 \times 7^3 \times 11$.
As $p\ge 17$, we obtain a contradiction.
\end{proof}

The space $S_2^{\mathrm{new}}(3^3 \cdot 3391)$ has dimension
$4520$. 
Using \texttt{Magma}, we compute the conjugacy classes
of eigenforms belonging
to this space and find that these have degrees
$1$, $1$, $1$, $1$, $554$, $556$, $564$, $564$, $565$, $565$, 
$574$ and $574$. The four rational eigenforms, of course, correspond
to the four elliptic curves $F_i$.
Unfortunately we have found it impossible to 
compute the coefficients of the irrational eigenforms due to the
enormous size of their fields of coefficients.
For a prime $\ell \ne 3$, $3391$, write $T_\ell$ for the Hecke operator
acting on $S_2^{\mathrm{new}}(3^3 \cdot 3391)$, and let
$C_\ell \in \Z[t]$ be the characteristic polynomial of $T_\ell$
(i.e. the $\ell$-th Hecke polynomial);
this is a polynomial of degree $4520$.
Using \texttt{Magma}, we found it straightforward (though
somewhat time-consuming) to compute
the polynomials $C_\ell(t)$ for $\ell <100$.
The polynomial $C_\ell$ satisfies 
\[
C_\ell(t)=\prod (t-a_\ell(f))
\]
where $f$ runs through the eigenforms of weight $2$ and level 
$3^3 \cdot 3391$. Note that $t-a_\ell(F_i)$ divides $C_\ell(t)$
for $i=1$, $2$, $3$ and $4$. We let
\[
C^\prime_\ell(t)=\frac{C_\ell(t)}{\prod_{1 \le i \le 4} (t-a_\ell(F_i))} \, .
\]
We now let
\[
R_\ell(a)=
\begin{cases}
C^\prime_\ell(\ell+1) \cdot C^\prime_\ell(\ell-1) & \text{if $\ell \mid \Delta_a$},\\
C^\prime_\ell(a_\ell(E_a)) & \text{if $\ell \nmid \Delta_a$}.
\end{cases}
\]
If $\ell \ne 2$, we let
\[
B_\ell=\ell \cdot \prod_{0 \le a \le \ell} R_\ell(a)
\]
and set $B_2=C_2^\prime(2)$.
\begin{lem}\label{lem:Bell}
Let $(x,z,p)$ be a solution to \eqref{eqn:k6} with $p \ge 17$
prime. Let $\ell \ne 3$, $3391$ be prime. Then $p \mid B_\ell$.
\end{lem}
\begin{proof}
By Lemmata~\ref{lem:ll} and~\ref{lem:notrat},
we know that $E_x \sim_p f$ where $f$ is an irrational eigenform
of weight $2$ and level $3^3 \cdot 3391$. It follows from
the above that $t-a_\ell(f)$ is a factor of $C_\ell^\prime$.
The lemma now follows from Lemma~\ref{lem:info} (for 
$\ell=2$ we are making use of the fact that $E_x$ has
good reduction at $2$ and that $a_2(E_x)=2$).
\end{proof}

\begin{proof}[Proof of Theorem~\ref{thm:2}]
Let $(x,z,p)$ be a solution to \eqref{eqn:k6} with $p \ge 11$.
Let 
$$
P=\{2\} \cup \{5,7,11,\dotsc,97\}
$$
be the set of primes less than $100$, excluding $3$.
Using \texttt{Magma},
we find that
\[
\gcd(\{B_\ell \; : \; \ell \in P\}) = 2^{27} \cdot 3^{28} \cdot 5^3 \cdot 7.
\]
This computation took roughly $21$ hours on a 2500MHz AMD Opteron.
The computation time was dominated by the computation of the polynomials
$C^\prime_\ell$.
The desired result then follows from Lemma~\ref{lem:Bell}.  
\end{proof}
We remark in passing
that the integers $B_\ell$ are extremely large which is why we
do not reproduce any of them here. By way of example, $\lvert B_2 \rvert \approx
1.1 \times 10^{569}$.

\section{The equation $(x-1)^k+x^k+(x+1)^k=y^p$ with $k \ge 7$}

It is natural to wonder if it is possible to attach 
a Frey-Hellegouarch curve a solution of the equation $(x-1)^k+x^k+(x+1)^k=z^p$
for exponents  $k \ge 7$. It is easy to see that
\[
(x-1)^k+x^k+(x+1)^k=\begin{cases}
f(x^2) & \text{if $k$ is even}\\
x f(x^2) & \text{if $k$ is odd},
\end{cases}
\]
where $f \in \Z[x]$. For $7 \le k \le 50$, say, we find that the polynomials
$f$ are irreducible and all their roots are real. We are unable to
prove that this is true in general for higher values of $k$. Suppose now
that $f$ is indeed a totally real irreducible polynomial, let $\theta$
be a root, and let $K=\Q(\theta)$. 
By a standard descent argument,
 $x^2-\theta=\alpha \xi^p$ where $\alpha$ belongs to a finite set
and $\xi$ is an integer in $K$. This can be viewed as a $(p,p,2)$
Fermat equation to which one can apply modularity and level-lowering
results over the totally real field $K$, in a similar manner to
that of several recent papers, e.g. 
\cite{BDMS}, \cite{DF}, \cite{FS}. We hope to pursue this approach
in a forthcoming paper.

\end{document}